\documentclass[11pt,twoside]{article}

\usepackage{amsmath}
\usepackage{amsthm}
\usepackage{amsfonts, amssymb}
\usepackage{mathrsfs}
\usepackage[all]{xy}

\def\titulo#1{\noindent{\bf\LARGE{#1}} \bigskip \thispagestyle{plain}}
\def\autor#1{\noindent{\sc #1}\smallskip}
\def\direccion#1{\noindent #1\bigskip}
\def\email#1{\vspace{1cm}\noindent E-mail address: {\sf #1} \bigskip}

\usepackage{latexsym}
\usepackage[dvips]{graphicx}

\theoremstyle{plain}
\newtheorem{lema}{Lemma}[section]
\newtheorem{prop}[lema]{Proposition}
\newtheorem{teo}[lema]{Theorem}

\newtheorem{fmteo}[lema]{First Main Theorem}
\newtheorem{smteo}[lema]{Second Main Theorem}
\newtheorem{coro}[lema]{Corollary}
\theoremstyle{remark}

\newtheorem{obs}[lema]{Remark}
\newtheorem{obsi}[lema]{Important Remark}
\newtheorem{notac}[lema]{Notation}
\theoremstyle{definition}
\newtheorem{defi}[lema]{Definition}
\newtheorem{ej}[lema]{Example}

\def\ce{\searrow \! \! \! \! \! ^e \: \: }
\def\ee{\nearrow \! \! \! \! \! \! \! ^e \ \ }
\def\he{\overset{he}{\simeq}}
\def\we{\overset{we}{\approx}}
\def\se{\diagup \! \! \! \! \: \! \: \searrow}

\def\e{\mathfrak{e}}

\def\k{\mathcal{K}}
\def\x{\mathcal{X}}
\def\c{\mathcal{C}}
\def\d{\mathcal{D}}
\def\dd{\overline{\mathcal{D}}}
\def\e{\mathcal{E}}
\def\etft0{\mathnormal{etfT_0}}

\begin{document}

\titulo{Simple Homotopy Types and Finite Spaces}

\autor{Jonathan Ariel Barmak, Elias Gabriel Minian}

\direccion{Departamento  de Matem\'atica.\\
 FCEyN, Universidad de Buenos Aires. \\ Buenos
Aires, Argentina}

\begin{abstract}
\noindent We present a new approach to simple homotopy theory of polyhedra using finite topological spaces. We define the concept of \it collapse \rm 
of a finite
space and prove that this new notion corresponds exactly to the concept of a simplicial collapse. 
More precisely, we show that a collapse $X\searrow Y$ of finite spaces induces a 
simplicial collapse $\k(X)\searrow \k(Y)$ of their associated simplicial complexes. Moreover, 
a simplicial collapse $K\searrow L$ induces a collapse $\x(K)\searrow \x(L)$ of the associated finite spaces. This establishes
 a one-to-one correspondence between simple homotopy types of finite simplicial complexes and simple equivalence classes of finite spaces.

\noindent We also prove a similar result for maps: We give a complete characterization of the class of maps between finite spaces 
which induce simple homotopy equivalences between the
associated polyhedra. Furthermore, this class describes all maps coming from simple homotopy equivalences at the level of complexes.

\noindent The advantage of this theory is that the elementary move of finite spaces is much simpler than the elementary 
move of simplicial complexes: It consists of
removing (or adding) just a single point of the space. 

\end{abstract}



\section{Introduction}

J.H.C. Whitehead's theory of simple homotopy types is inspired by Tietze's theorem in combinatorial group theory, which states that any finite presentation
of a group could be deformed into any other by a finite sequence of elementary moves, which are now called Tietze transformations. Whitehead 
translated these algebraic moves into the well-known geometric moves of elementary collapses and expansions of finite simplicial complexes. 
His beautiful theory of simple homotopy types turned out to be fundamental for the development of piecewise-linear topology: 
The s-cobordism theorem, 
Zeeman's conjecture \cite{Zee}, the applications of the theory in surgery, Milnor's classical paper on Whitehead Torsion \cite{Mil} and the topological 
invariance of torsion are some of its major uses and advances.

\medskip

In this paper we show how to use finite topological spaces to study simple homotopy types. There is a strong relationship between finite spaces and
finite simplicial complexes, which was discovered by McCord \cite{Mcc}. Explicitly, given a finite simplicial complex $K$, one can associate to $K$ 
a finite $T_0$-space $\x (K)$ which corresponds to the poset of simplices of $K$ ordered by inclusion. Moreover, a simplicial map $\varphi:K\to L$ gives rise to 
a continuous map $\x(\varphi)$ between the associated finite spaces. Conversely, one can associate to a finite $T_0$-space $X$ a simplicial complex $\k (X)$, 
whose simplices are the non-empty chains of $X$, and a weak equivalence $\k (X)\to X$. This construction is also functorial.

\medskip

In \cite{Bar} we showed that finite spaces are very useful for studying homotopy invariants of (general) spaces. In fact, in that
article we were looking for \it finite minimal models \rm of some spaces, i.e. the smallest finite spaces which are weak equivalent to a given space.  
In our opinion, finite spaces are more effective for studying homotopy theory  than simplicial complexes, because, besides their combinatorial 
nature, they have the extra topological structure.

\medskip

It is easy to prove that if two finite $T_0$-spaces $X, Y$ are homotopy equivalent, their associated simplicial complexes $\k(X), \k(Y)$ are 
also homotopy equivalent. Furthermore, Osaki \cite{Osa} showed that in this case, the latter have the same simple homotopy type. 
Nevertheless, we noticed that the converse of this result is not true in general: There are finite spaces with different homotopy types 
whose associated simplicial complexes have the same simple homotopy type. Starting from this point, 
we were looking for the relation that $X$ and $Y$ should satisfy for their associated complexes to be simple homotopy equivalent.
More especifically, we wanted to find an elementary move  in the setting of finite spaces (if it existed) which corresponds exactly to 
 a simplicial collapse of the associated polyhedra.

\medskip

We discovered this elementary move when we were looking for a homotopically trivial finite space (i.e. weak equivalent to a point) which were non-contractible.
In order to construct such a space, we developed a method of reduction (i.e. a method that allows us to reduce a finite space to a smaller weak equivalent
space). This method of reduction together with the homotopically trivial and non-contractible space (of 11 points) that we found are exhibited in section
3. Suprisingly, this method, which consists of removing a \it weak point \rm  of the space (see \ref{weakpoint}), turned out to be the key to solve 
the problem of translating 
simplicial collapses into this setting.

\medskip

We will say that two finite spaces are \textit{simply equivalent} if we can obtain one of them  from the other by adding and removing weak points. 
If $Y$ is obtained from $X$ by only removing weak points, we say that $X$ \textit{collapses} to $Y$ and write $X \searrow Y$.
The first main result of this article is the following

\begin{fmteo}

\begin{enumerate}

\item[ ]
\item[(a)] Let $X$ and $Y$ be finite $T_0$-spaces. Then, $X$ and $Y$ are simply equivalent if and only if $\k(X)$ and $\k (Y)$ have the same simple homotopy type. 
Moreover, if $X \searrow Y$ then $\k (X) \searrow \k(Y)$.
\item[(b)] Let $K$ and $L$ be finite simplicial complexes. Then, $K$ and $L$ are simple homotopy equivalent if and only if $\x(K)$ and $\x (L)$ are simply equivalent. Moreover, if $K \searrow L$ then $\x (K) \searrow \x (L)$.
\end{enumerate}
\end{fmteo}

In particular, the functors $\k$ and $\x$ induce a one-to-one correspondence between simply equivalence classes of finite spaces and 
simple homotopy types:

\begin{displaymath}
\xymatrix@C=50pt{ \{Finite\ T_0-Spaces\} \! \! \textrm{\raisebox{-2ex}{\Huge{/}} \raisebox{-2.7ex}{$\! \! \! \! \! \se$}}
\ar@<2.4ex>^{\! \! \! \! \! \! \! \! \! \! \! \! \k}[r] &  \{Finite\ Simplicial\ Complexes\} \! \! \textrm{\raisebox{-2ex}{\Huge{/}} \raisebox{-2.7ex}{$\! \! \! \! \! \se$}} \ar@<-0.3ex>[l]^{\! \! \! \! \! \! \! \! \! \! \! \! \x} } 
\end{displaymath}

We are now able to study finite spaces using all the power of Whitehead's simple homotopy theory for CW-complexes. But also, what is more important, 
we can use finite spaces to strengthen the classical theory. 
The elementary move in this setting  is much simpler to handle and describe because it consists of adding or removing just one single point.

\medskip

As an example or application  of this  theorem,  we study  \textit{collapsible} finite spaces and their relationship with collapsible complexes. 
We also relate simple types of finite spaces with the notion of minimal finite model introduced in \cite{Bar}.

\medskip

In the last section of this article we investigate the class of maps between finite spaces which induce simple homotopy equivalences between their associated simplicial
 complexes. To this end, we introduce the notion of a \it distinguished \rm map. Similarly to the classical case, the class of maps which induce simple 
homotopy equivalences can be generated, in a certain way, by expansions and a kind of formal homotopy inverses of expansions. Remarkably this class, which we denote $\dd$,  
is also generated by the
 distinguished maps. The second main result of the article is the following
 
\begin{smteo}
\begin{enumerate}
\item[]
\item[(a)] Let $f:X\to Y$ be a map between finite $T_0$-spaces. Then $f\in \dd$ if and only if $\k(f):\k(X)\to \k(Y)$ is a simple homotopy equivalence. 

\item[(b)] Let $\varphi:K\to L$ be a simplicial map between finite simplicial complexes. Then $\varphi$ is a simple homotopy equivalence if and only if $\x(\varphi)\in \dd$.
\end{enumerate}
\end{smteo} 
  
\section{Preliminaries}

In this section we recall  the basic notions on finite spaces which are essential in section 3 and 4. For more details on finite spaces 
we refer the reader to \cite{Mcc, Sto} and P. May's beautiful notes \cite{May, May2}.

First we describe Alexandroff's correspondence between topologies and preorders on a finite set.

\bigskip

Let $X$ be a finite topological space and let $x$ be a point of $X$. The \textit{minimal open set} $U_x$ of $x$ is defined as the intersection of all open sets containing $x$. The minimal open sets of $X$ form a basis for the topology on $X$, which is called the \textit{minimal basis} of $X$ for obvious reasons. 

The preorder associated to the topology on $X$ is given by the relation $x\le y$ if $x\in U_y$.

\bigskip

Conversely, given a preorder $\le$ on $X$, we define for each $x\in X$ the set $$U_x=\{y\in X \ | \ y\le x\}.$$
It is not hard to see that these sets form a basis for a topology, which is the topology associated to $\le$.

\bigskip

These applications define a one-to-one 
correspondence between topologies and preorders in the finite set $X$. Moreover, $T_0$-topologies correspond to orders. One can therefore regard
finite $T_0$-spaces as finite posets and viceversa.

\bigskip

It is very useful to  represent finite $T_0$-spaces with their \textit{Hasse diagram}. The Hasse diagram of $X$ is a digraph whose vertex set is $X$ and whose edges are the ordered pairs $(x,y)$ such that $x<y$ and there exist no $z\in X$ with $x<z<y$.

\bigskip

Consider for example the space $X=\{a,b,c,d\}$ whose proper open sets are $\{a,c,d\}$, $\{b,c,d\}$, $\{c,d\}$ and $\{d\}$. Its Hasse diagram would be

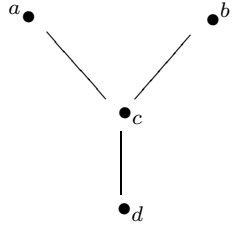
\begin{displaymath}
\xymatrix@C=10pt{ ^a \bullet \ \ \ar@{-}[dr] & & \ \ \bullet ^b \ar@{-}[ld] \\
									& \ \ \bullet _c \ar@{-}[d] &  \\
									& \ \ \bullet _d \\ } 
\end{displaymath}

Instead of representing an edge $(x,y)$ with an arrow, one simply writes $y$ over $x$.

\bigskip

Given a finite space $X$, we will denote by $X^{op}$ the space with the same underlying set as $X$ but with the opposite preorder. 

\bigskip

\bigskip

Following Stong \cite{Sto} and May \cite{May},  we say that a point $x$ of a finite $T_0$-space $X$ is an \textit{up beat point} if the set of points which 
are greater than $x$ has a minimum. A \textit{down beat point} is one such that the set of points below it has a maximum.\

\begin{obs} \label{comparable}
If $x\in X$ is a beat point (up or down), there exists $y\in X$, $y\neq x$ such that any point which is comparable with $x$ is also comparable with $y$.
\end{obs}

Stong proved in \cite{Sto} that if $x\in X$ is a beat point, then $X\smallsetminus \{x\}$ is a strong deformation retract of $X$.\

A finite $T_0$-space with no beat points is called a \textit{minimal finite space}. A \textit{core} of a finite space $X$ is a strong deformation 
retract of $X$ which is a minimal finite space. \

Given a finite $T_0$-space $X$, one can find a sequence of spaces $X=X_0\supsetneq X_1 \supsetneq \ldots \supsetneq X_n$ 
where $X_{i+1}$ is obtained from $X_i$ by removing a beat point and such that $X_n$ has no beat points.
Therefore every finite $T_0$-space has a core. 

\smallskip

Furthermore, Stong proved that 
every homotopy equivalence between minimal finite spaces is a homeomorphism and therefore, the core of any finite space $X$ is unique up to homeomorphism. It can be
described as the smallest space which is homotopy equivalent to $X$.

\begin{obs}
If $X$ is a contractible finite $T_0$-space, there exists a sequence $X=X_0\supsetneq X_1\supsetneq \ldots \supsetneq X_n=*$ where $X_{i+1}$ is obtained from $X_i$ by removing a beat point.
\end{obs}

\smallskip

It is not hard to prove that a finite $T_0$-space with maximum or minimum is contractible.

\begin{obs}
A point $x$ of a finite $T_0$-space $X$ is an up beat point if and only if $x$ is a down beat point of $X^{op}$. Therefore, $X$ is contractible if and only if $X^{op}$ is contractible.
\end{obs}

\bigskip

Note that a function $f:X \to Y$ between finite spaces is continuous if and only if it is order preserving.\

Given two functions $f,g:X \to Y$, we will say that $f\le g$ if $f(x)\le g(x)$ for every $x\in X$.  It is not difficult to prove that 
if $f$ and $g$ are continuous and $f\le g$, then $f$ is homotopic to $g$ (see \cite{May,Sto} for more details).

\bigskip

\bigskip

Following McCord \cite{Mcc} (cf. also \cite{Bar, May2}) one can associate to a finite $T_0$-space $X$ a simplicial complex $\k (X)$, whose simplices are the non-empty chains of $X$, and a weak equivalence $|\k (X)|\to X$. Here $|\k(X)|$ denotes the geometric realization of $\k(X)$.

The application $\k$ is in fact functorial. A continuous map $f:X\to Y$ between finite $T_0$-spaces induces a simplicial map $\k(f):\k (X)\to \k(Y)$ which coincides with $f$ on vertices.

\bigskip

Given a map $f:X\to Y$, one has the following commutative diagram

\begin{displaymath}
\xymatrix@C=20pt{ |\k (X)| \ar@{->}[d] \ar@{->}^{|\k(f)|}[r] & |\k (Y)| \ar@{->}[d] \\
									X \ar@{->}^f[r] & Y } 
\end{displaymath}

If $f\simeq g: X\to Y$, it can be proved that $\k(f),\k(g) :\k (X)\to \k(Y)$ lie in the same contiguity class. In particular $|\k(f)|\simeq |\k(g)|$.

\bigskip 

Conversely, one can associate to each finite simplicial complex $K$ a finite $T_0$-space $\x (K)$. This finite space is the poset of simplices of $K$ ordered by inclusion. 

Note that $\k (\x (K))=K'$ is the first barycentric subdivision of $K$. This implies of course that there exists a weak equivalence $|K|\to \x (K)$.

This application is also functorial. A simplicial map $\varphi:K\to L$ between finite simplicial complexes induces a continuous map $\x(\varphi): \x (K)\to \x (L)$, where $\x(\varphi)(S)=\varphi(S)$ for every simplex $S$ of $K$.

In this case one has the following diagram that commutes up to homotopy

\begin{displaymath}
\xymatrix@C=20pt{ |K| \ar@{->}[d] \ar@{->}^{|\varphi|}[r] & |L| \ar@{->}[d] \\
									\x (K) \ar@{->}^{\x(\varphi)}[r] & \x (L) } 
\end{displaymath}

\bigskip

Recall that two spaces $X$ and $Y$ (non-necessarily finite) are weak equivalent if there exists a sequence of spaces $X=X_1, X_2, \ldots ,X_n=Y$ such that for each $1\le i<n$ there is a weak equivalence $X_i\to X_{i+1}$ or $X_{i+1}\to X_i$. We denote $X\we Y$.

\bigskip

We call a space $X$ homotopically trivial if $X\we *$ (i.e. all its homotopy groups are trivial).

\bigskip 

Note that if $X$,$Y$ are finite $T_0$-spaces, then by Whitehead Theorem,  $X\overset{we}{\approx} Y$ if and only if $|\k (X)|$ and $|\k(Y)|$ have the same homotopy type.

For finite simplicial complexes $K$ and $L$, $|K|$ and $|L|$ are homotopy equivalent if and only if $\x (K)\overset{we}{\approx} \x (L)$.

\bigskip

\bigskip

We finish this introductory section by recalling  the basic notions on  simple homotopy theory for simplicial complexes. Essentially we want to fix 
the notations that we will use in the main sections of the paper. The standard references for this are of course Whitehead's papers
\cite{Whi, Whi2, Whi4}, Milnor's article \cite{Mil} and M.M.Cohen's book \cite{Coh}.

\bigskip

Let $K$ and $L$ be finite simplicial complexes. Recall that there is an \textit{elementary simplicial collapse} from $K$ to $L$ if there is a simplex 
$S$ of $K$ and a vertex $a$ of $K$ not in $S$ such that $K=L\cup aS$ and $L\cap aS=a\dot{S}$. 
Elementary collapses will be denoted, as usual, $K\ce L$.

\bigskip

We say that $K$ \textit{(simplicially) collapses} to $L$ (or that $L$ \textit{expands} to $K$) if there exists 
a sequence $K=K_1, K_2, \ldots, K_n=L$ of finite simplicial complexes such that $K_i\ce K_{i+1}$ for all $i$. This is denoted by $K \searrow L$ or $L \nearrow K$.  

Two complexes $K$ and $L$ have the same simple homotopy type if there is a sequence $K=K_1, K_2, \ldots, K_n=L$ such that $K_i\searrow K_{i+1}$ 
or $K_i \nearrow K_{i+1}$ for all $i$. Following M.M. Cohen's notation, we denote this by $K \se L$.

\bigskip

It is well known that $K \se L$ if and only if $|K|$ and $|L|$ are simple homotopy equivalent  viewed as CW-complexes \cite{Whi4}.

\section{Simple Homotopy Types: The First Main Theorem}

We start by introducing  the notion of a \textit{weak beat point}. This concept appeared naturally when we were searching for reduction methods (cf. 
\cite{Bar}) to find the \it minimal finite models  \rm of some spaces, i.e. the smallest finite spaces which are weak equivalent to a given space.
Surprisingly this new notion turned out to be crucial for studying the relationship between finite spaces and simplicial collapses of their associated complexes.

\begin{defi} \label{weakpoint}
Let $X$ be a finite $T_0$-space. We will say that $x\in X$ is a \textit{weak beat point of $X$} (or a \textit{weak point}, for short) 
if either $U_x\smallsetminus \{x\}$ is contractible or $\overline{\{x\}} \smallsetminus \{x\}$ is contractible.
\end{defi}

Here $\overline{\{x\}}$ denotes the closure of $\{x\}$, i.e. the set of points which are greater than or equal to $x$.
Note that if $x\in X$ is a down beat point, $U_x\smallsetminus \{x\}$ has a maximum and if $x$ is an up beat point, $\overline{\{x\}} \smallsetminus \{x\}$ has a minimum. Therefore, beat points are particular cases of weak points.

\bigskip

As we have seen in the previous section, when $x$ is a beat point of $X$, the inclusion $i: X\smallsetminus \{x\} \hookrightarrow X$ is a homotopy equivalence. The following proposition generalizes Stong's result.

\begin{prop} \label{weak}
Let $x$ be a weak point of a finite $T_0$-space $X$. Then the inclusion map $i: X\smallsetminus \{x\} \hookrightarrow X$ is a weak equivalence.
\end{prop}

\begin{proof}
We may suppose that $U_x\smallsetminus \{x\}$ is contractible since the other case follows from this one considering $X^{op}$. Note that the minimal open set $U_x$ of $x$ in $X^{op}$ is the closure $\overline{\{x\}}$ of $x$ in $X$ and that $\k(X^{op})=\k(X)$. 

Given $y\in X$, we have that $i^{-1}(U_y)=U_y \smallsetminus \{x\}$, which has maximum if $y\neq x$ and is contractible if $y=x$. It is clear 
then that $$i|_{i^{-1}(U_y)}:i^{-1}(U_y)\to U_y$$ is a weak homotopy equivalence for every $y\in X$. Now the result follows from Theorem 6 
of \cite{Mcc} applied to the basis-like cover given by the minimal basis of $X$.
\end{proof}

In \cite{Osa} Corollary 3.4, Osaki proves that if two finite $T_0$-spaces, $X$ and $Y$ are homotopy equivalent, then their associated simplicial complexes, $\k (X)$ and $\k (Y)$, have the same simple homotopy type.

\bigskip

One might ask if the converse of this result also holds. Explicitly, suppose that $\k (X)$ and $\k(Y)$ have the same simple 
homotopy type, is it true that $X$ and $Y$ are homotopy equivalent?

This question, which we can refer to as (Q1), is related to the following question (Q2): Is there a finite space 
which is homotopically trivial but not contractible?

In \cite{Bar} we have already showed that Whitehead's Theorem does not hold for 
finite spaces (in fact, one can exhibit many examples of finite spaces which are weak equivalent but not homotopy equivalent), but 
we did not know whether (Q2) was true or not.

Both questions are related in the following sense: An affirmative answer to (Q2) would give a negative answer to (Q1). If $X$ is a homotopically trivial finite $T_0$-space, then $|\k (X)|$ is contractible, then its Whitehead group is trivial, and therefore, $\k (X) \se *$. 

\smallskip

As an application of the last proposition, we found the following example of a homotopically trivial space of 11 points which is not contractible.

\begin{ej}[The Wallet] \label{wallet}

Let us consider $W$, the finite $T_0$-space whose Hasse diagram is the following

\begin{displaymath}
\xymatrix@C=10pt{ \bullet \ar@{-}[d] \ar@{-}[drr] & & \bullet \ar@{-}[lld] \ar@{-}[rrd] & & \underset{}{\overset{x}{\bullet}} \ar@{-}[lld] \ar@{-}[rrd] & & \bullet \ar@{-}[lld] \ar@{-}[d]  \\
		\bullet \ar@{-}[dr] \ar@{-}[drrr] & & \bullet \ar@{-}[dl] \ar@{-}[dr] & & \bullet \ar@{-}[dl] \ar@{-}[dr] & & \bullet \ar@{-}[dlll] \ar@{-}[dl] \\
		& \bullet & & \bullet & & \bullet } 
\end{displaymath}

\begin{center}
Fig. 1: $W$
\end{center}

Note that $W$ is not contractible since it is a minimal finite space (with more than one point). Nevertheless it contains a weak point $x$ (see Fig. 1), since $U_x \smallsetminus \{x\}$ is contractible (see Fig. 2).

\begin{displaymath}
\xymatrix@C=10pt{ & \bullet \ar@{-}[dl] \ar@{-}[dr] & & & & \bullet \ar@{-}[dlll] \ar@{-}[dl] \\
		 \bullet & & \bullet & & \bullet } 
\end{displaymath}

\begin{center}
Fig. 2: $U_x\smallsetminus \{x\}$
\end{center}

Therefore $W$ is weak equivalent to $W\smallsetminus \{x\}$, whose Hasse diagram is as follows

\begin{displaymath}
\xymatrix@C=10pt{ \bullet \ar@{-}[d] \ar@{-}[drr] & & \bullet \ar@{-}[lld] \ar@{-}[rrd] & & & & \bullet \ar@{-}[lld] \ar@{-}[d]  \\
		\bullet \ar@{-}[dr] \ar@{-}[drrr] & & \bullet \ar@{-}[dl] \ar@{-}[dr] & & \bullet \ar@{-}[dl] \ar@{-}[dr] & & \bullet \ar@{-}[dlll] \ar@{-}[dl] \\
		& \bullet & & \bullet & & \bullet } 
\end{displaymath}

\begin{center}
Fig. 3: $W\smallsetminus \{x\}$
\end{center}

Now it is easy to see that this subspace is contractible. In fact, $W \smallsetminus \{x\}$ does have beat points, and one can get rid of them one by one.

Therefore $W$ is homotopically trivial but not contractible.

\end{ej}

As we pointed out before, this example shows that the converse of Osaki's result does not hold.

This naturally leads us to the following harder problem, which is one of the main goals of the paper: What is exactly the 
relation that $X$ and $Y$ must satisfy for $\k(X)$ and $\k (Y)$ to have the same simple homotopy type?

More precisely, is there an \textit{elementary move} in the setting of finite spaces which corresponds to a 
simplicial collapse of the associated complexes?

\smallskip

We found out that the notion of weak point (\ref{weakpoint}) was the key  to solve this problem:

\begin{defi} \label{definicion}
Let $X$ be a finite $T_0$-space and let $Y\subsetneq X$. We say that $X$ \textit{collapses} to $Y$ by an \textit{elementary collapse} 
(or that $Y$ \textit{expands} to $X$ by an \textit{elementary expansion}) if $Y$ is obtained from $X$ by removing a weak point. We denote $X\ce Y$ or $Y\ee X$.

\smallskip

In general, given two finite $T_0$-spaces $X$ and $Y$, we say that $X$ \textit{collapses} to $Y$ (or $Y$ \textit{expands} to $X$) if there is a sequence $X=X_1, X_2, \ldots, X_n=Y$ of finite $T_0$-spaces such that for each $1\le i <n$, $X_i\ce X_{i+1}$. In this case we write $X\searrow Y$ or $Y\nearrow X$.

\smallskip

Two finite $T_0$-spaces $X$ and $Y$ are \textit{simply equivalent} if there is a sequence $$X=X_1, X_2, \ldots, X_n=Y$$
 of finite $T_0$-spaces such that for each $1\le i <n$, $X_i\searrow X_{i+1}$ or $X_i\nearrow X_{i+1}$. We denote in this case $X \se Y$, 
following the same notation that we adopted for simplicial complexes.
\end{defi}

It follows straightforward from \ref{weak} that, if $X$ and $Y$ are simply equivalent finite $T_0$-spaces, then they are weak equivalent.

We shall see later, as an immediate corollary of our first main result \ref{main}, that homotopy equivalent finite spaces are simply equivalent. This result follows from the fact that beat points are weak points and that homeomorphic finite spaces are simply equivalent.

The main reason of this curious situation (in the classical setting, a simple homotopy equivalence is in particular a homotopy equivalence) is that 
Whitehead's Theorem does not hold in this context.

\bigskip

In order to prove the First Main Theorem, we need some previous results. We show first that the associated finite space $\x (K)$ of a simplicial cone $K$ is contractible.

\bigskip

Suppose $K=aL$ is a cone, i.e. $K$ is the join of a simplicial complex $L$ with a vertex $a\notin L$. Since $|K|$ is contractible, it is clear that $\x (K)$ is homotopically trivial. The following lemma shows that $\x (K)$ is, in fact, contractible.

\begin{lema} \label{cono}
Let $K=aL$ be a finite cone. Then $\x (K)$ is contractible.
\end{lema}
\begin{proof}
Define $f:\x (K)\to \x (K)$ by $f(S)=S\cup \{a\}$. This function is order-preserving and therefore continuous. 

If we consider $g:\x (K)\to \x (K)$ the constant map that takes all $\x (K)$ into $\{a\}$, we have that $$ 1_{\x (K)}\le f \ge g .$$
This proves that the identity is homotopic to a constant map.

\end{proof}

It is well known that any finite simplicial complex $K$ has the same simple homotopy type of its barycentric subdivision $K'$. We prove next an analogous result for finite spaces.

\bigskip

Following \cite{Har}, the barycentric subdivision of a finite $T_0$-space $X$ is defined by $X'=\x (\k (X))$. Explicitly, $X'$ consists
 of the non-empty chains of $X$ ordered by inclusion. It is shown in \cite{Har} that there is a weak equivalence $X'\to X$ which takes each chain $C$ to $max(C)$.

\bigskip

\begin{notac}
We fix some notation that we will adopt for the rest of the paper.

Given a set $A$ which can be regarded as a subset of different spaces $X$ and $Y$, we will denote $\overline{A}^X$ the closure of $A$ in the space $X$ 
so as not to confuse it with its closure in $Y$. Similarly, given $x\in X$, we will denote $U_x^X$ the minimal open set of $x$ in the space $X$.
\end{notac}

\begin{prop} \label{bedex}
Let $X$ be a finite $T_0$-space. Then $X$ and $X'$ are simply equivalent.
\end{prop}

\begin{proof}
Consider the space $B(X)$ whose underlying set is $X\sqcup X'$ with the following order relation. Given $a,b\in B(X)$ we say that $a\le b$ if one of the following holds:

\begin{itemize}
\item $a,b\in X$ and $a\le b$ in $X$.
\item $a,b\in X'$ and $a\le b$ in $X'$.
\item $a\in X'$, $b\in X$ and $max(a)\le b$ in $X$. 
\end{itemize}

It easy to see that this relation defines an order on $B(X)$, thus it is a finite $T_0$-space.

We will show that $X$ and $X'$ are both simply equivalent to $B(X)$.

\bigskip

We label all the elements $C_1, C_2, \ldots, C_n$ of $X'$ in such a way that $C_i\le C_j$ implies $i\le j$. Then we define $X_i=X\sqcup \{C_1, C_2, \ldots, C_i\} \subseteq B(X)$ for every $0\le i\le n$. 

Since $$\overline{\{C_i\}}^{B(X)}\smallsetminus \{C_i\}=\{x\in X \ | \ x\ge max(C_i)\}\sqcup \{C\in X' \ | \ C\supsetneq C_i\}, $$

we have that $$\overline{\{C_i\}}^{X_i}\smallsetminus \{C_i\}=\{x\in X \ | \ x\ge max(C_i)\},$$

which is homeomorphic to $\overline{\{max(C_i)\}}^X$. Therefore, it has a minimum and then is contractible.

We have just proved that $C_i$ is a weak point of $X_i$ for every $1\le i\le n$. Hence, $X_i \ce X_{i-1}$ for $1\le i\le n$, and then $X=X_0$ is simply equivalent to $B(X)=X_n$.

\bigskip

Now order the elements $x_1,x_2, \ldots, x_m$ of $X$ in such a way that $x_i\le x_j$ implies $i\le j$. We define $X_i'=\{x_{i+1},x_{i+2}, \ldots, x_m\}\sqcup X'\subseteq B(X)$ for every $0\le i\le m$. Then $$U_{x_{i}}^{B(X)}\smallsetminus \{x_{i}\}=\{x\in X \ | \ x< x_{i}\}\sqcup \{C\in X' \ | \ max(C)\le x_{i}\}$$

Therefore $$U_{x_{i}}^{X_{i-1}'}\smallsetminus \{x_{i}\}=\{C\in X' \ | \ max(C)\le x_{i}\},$$ which is homeomorphic to $\x (\k (U_{x_{i}}^X))$.

But $\k (U_{x_{i}}^X)=x_{i}\k (U_{x_{i}}^X\smallsetminus \{x_{i}\})$ is a cone. By the previous lemma, $U_{x_{i}}^{X_{i-1}'}\smallsetminus \{x_{i}\}$ is contractible. Thus $x_{i}$ is a weak point of $X_{i-1}'$ for every $1\le i\le m$, and then $X_{i-1}'\ce X_i'$ for $1\le i\le m$. Therefore $B(X)=X_0'$ is simply equivalent to $X'=X_m'$.
\end{proof}

The next technical lemma will also be used in the proof of the First Main Theorem.

\begin{lema} \label{expansion}
Let $L$ be a subcomplex of a finite simplicial complex $K$. Let $T$ be a set of simplices of $K$ which are not in $L$, 
and let $a$ be a vertex of $K$ which is contained in no simplex of $T$, but such that $aS$ is a simplex of $K$ for every $S\in T$. 
Finally, suppose that $K=L\cup \bigcup\limits_{S\in T} \{S,aS\}$ (i.e. the simplices of $K$ are those of $L$ together with the simplices $S$ and $aS$ for  
every $S$ in
$T$). Then $L$  simplicially expands  to $K$. 
\end{lema}
\begin{proof}
Suppose that $T=\{S_1,S_2, \ldots, S_n\}$ where $i\le j$ implies $\# S_i\le \# S_j$.

We define $K_i=L\cup \bigcup\limits_{j=1}^{i} \{S_j,aS_j\}$ for $0\le i \le n$. 

Let $1\le i\le n$, and let $S\subsetneq S_i$. If $S\in T$, since $\# S<\# S_i$, we have that $S, aS\in K_{i-1}$. If $S\notin T$, then $S, aS\in L\subseteq K_{i-1}$. Therefore, we have proved that $aS_i\cap K_{i-1}=a\dot{S_i}$. 

\smallskip

Inductively we have that $K_i$ is a simplicial complex for every $i$ and that there is an elementary simplicial expansion from $K_{i-1}$ to $K_i$ for every $1\le i\le n$, thus $L=K_0$ expands simplicially to $K=K_n$.  
\end{proof}

Now we are ready to prove the first main result of this article.

\begin{fmteo} \label{main}
\begin{enumerate} 
\item[ ]
\item[(a)] Let $X$ and $Y$ be finite $T_0$-spaces. Then, $X$ and $Y$ are simply equivalent if and only if $\k(X)$ and $\k (Y)$ have the same simple homotopy type. Moreover, if $X \searrow Y$ then $\k (X) \searrow \k(Y)$.
\item[(b)] Let $K$ and $L$ be finite simplicial complexes. Then, $K$ and $L$ are simple homotopy equivalent if and only if $\x(K)$ and $\x (L)$ are simply equivalent. Moreover, if $K \searrow L$ then $\x (K) \searrow \x (L)$.
\end{enumerate}
\end{fmteo}
\begin{proof}
Let $X$ be a finite $T_0$-space and let $x\in X$ be a weak point. Suppose first that $U_x\smallsetminus \{x\}$ is contractible. In this case, there exists a sequence of spaces  $U_x\smallsetminus \{x\}=X_n\supsetneq X_{n-1} \supsetneq \ldots \supsetneq X_1=\{x_1\}$, $X_i=\{x_1,x_2,\ldots ,x_i\}$ and such that $x_i$ is a beat point of $X_i$ for $i\ge 2$. 

It follows from \ref{comparable} that for each $2\le i\le n$ there exists $y_i\in X_{i-1}$ with the following property: if $z\in X_i$ is comparable with $x_i$, then it is comparable with $y_i$.

\bigskip

For every $1\le i \le n$ we define $K_i\subseteq \k (X)$, the subcomplex whose simplices are the chains of $Y=X\smallsetminus {\{x\}}$ and the chains of $\overline{\{x\}}\cup X_i\subseteq X$. In other words, $K_i=\k (Y) \cup \k (\overline{\{x\}}\cup X_i)$.

\bigskip

The simplices of $K_1$ which are not in $\k (Y)$, are the chains of $\overline{\{x\}}\cup \{x_1\}$ that 
contain $x$. Taking $T=\{S\in K_1 \ | \ x\in S, \ x_1\notin S\}$ and $a=x_1$ in \ref{expansion} it is easy to see 
that $\k (Y) \nearrow K_1$. If $S\in T$, every element of $S$ is greater than or equal to $x$, and therefore comparable with $x_1$. That is, $x_1S$ is a simplex of $K_1$.

\bigskip

Now, if $i\ge 2$, the simplices of $K_i$ which are not in $K_{i-1}$ are the chains of $\overline{\{x\}}\cup X_i$ that contain both $x$ and $x_i$. 

In order to use \ref{expansion}, we define $T=\{S\in K_i \ | \ x,x_i\in S, \ y_i\notin S\}$ and $a=y_i$. If $S\in T$, every element $y\in S$ satisfies one of the following: (i) $y\ge x$ or (ii) $y\in X_i$ is comparable with $x_i$. In any of these cases it holds that $y$ is comparable with $y_i$, and then $y_iS\in K_i$.

We deduce then that $K_{i-1}\nearrow K_i$ and therefore $\k(Y) \nearrow K_n=\k (X)$.

\bigskip 

Suppose now that $\overline{\{x\}}\smallsetminus \{x\}$ is contractible. If we consider the opposite order on X, it follows from the previous case that $\k (X\smallsetminus \{x\}) \nearrow \k(X)$.

We have then proved that $X\searrow Y$ implies $\k (X)\searrow \k(Y)$. In particular, $X\se Y$ implies $\k (X)\se \k(Y)$.

\bigskip

Now suppose that $K$ and $L$ are finite simplicial complexes such that there is an elementary simplicial collapse from $K$ to $L$. Hence, there exists $S\in K$ and a vertex $a$ of $K$ not in $S$ such that $aS\in K$, $K=L\cup \{S,aS\}$ and $aS\cap L=a\dot{S}$.

There is only one simplex of $K$ containing $S$ properly, namely $aS$. Therefore, $S$ is an up beat point in $\x (K)$, and then $\x (K) \ce \x (K)\smallsetminus \{S\}$. 

The simplices contained properly in $aS$ which are distinct from $S$ form a cone, and then $$U_{aS}^{\x (K)\smallsetminus \{S\}} \smallsetminus \{aS\}=\x (a\dot{S})$$ is contractible by \ref{cono}. Then $aS$ is a weak point of $\x (K)\smallsetminus \{S\}$ which collapses to $\x (K)\smallsetminus \{S,aS\}=\x (L)$.

This proves the first part of $(b)$ and the ``moreover'' part.

\bigskip

If $X$, $Y$ are finite $T_0$-spaces such that $\k (X)\se \k (Y)$, we have just proved that $$X'=\x (\k (X))\se Y'=\x (\k (Y)).$$
 However, by \ref{bedex} $X\se X'$ and $Y\se Y'$, and then $X\se Y$.\

Finally, if $K$, $L$ are finite simplicial complexes such that $\x (K) \se \x (L)$, then $$K'=\k (\x (K))\se L'=\k (\x (L)).$$ 
Since $K\se K'$ and $L\se L'$, it follows that $K\se L$. This completes the proof.
\end{proof}

In particular, we have the following corollary.

\begin{coro}
The functors $\k$, $\x$ induce a one-to-one correspondence between simply equivalence classes of finite spaces and simple homotopy types of finite simplicial complexes

\begin{displaymath}
\xymatrix@C=50pt{ \{Finite\ T_0-Spaces\} \! \! \textrm{\raisebox{-2ex}{\Huge{/}} \raisebox{-2.7ex}{$\! \! \! \! \! \se$}}
\ar@<2.4ex>^{\! \! \! \! \! \! \! \! \! \! \! \! \k}[r] &  \{Finite\ Simplicial\ Complexes\} \! \! \textrm{\raisebox{-2ex}{\Huge{/}} \raisebox{-2.7ex}{$\! \! \! \! \! \se$}} \ar@<-0.3ex>[l]^{\! \! \! \! \! \! \! \! \! \! \! \! \x} } 
\end{displaymath}
\end{coro}

The theorem shows that the exact translation of the notion of simple homotopy type in the context of finite spaces is the one defined in \ref{definicion}.

This result can be used in both ways. On one hand, to prove results on finite spaces using the machinery of Whitehead's  theory. On the other hand, to enrich the classical theory with the incipient theory of finite spaces.

We hope to find in the future a new way to characterize the obstruction for two simplicial complexes to have the same simple homotopy type, which is measured by the Whitehead groups of the complexes, via their associated finite spaces. Note that in the 
finite space setting the elementary move consists of removing (or adding) just one point.

\bigskip

If two finite $T_0$-spaces $X$ and $Y$ are homeomorphic, we could use a \mbox{construction similar} to $B(X)$ in \ref{bedex} to prove that $X\se Y$ (see \ref{dis}). However that is not necessary, \mbox{since this} \mbox{follows immediately} from the theorem because $\k (X)$ and $\k (Y)$ are isomorphic and therefore simple 
homotopy equivalent.

It follows that homotopy equivalent finite spaces are simply equivalent. Explicitly, if $X$ and $Y$ have the same homotopy type, their cores $X_c$ and $Y_c$ are homeomorphic and then $X\searrow X_c\se Y_c\nearrow Y$.

\bigskip 

The following diagrams illustrate the whole situation.
\begin{displaymath}
\xymatrix@C=20pt{ X \he Y  \ar@{=>}[r] & X \se Y \ar@{=>}[r] \ar@{<=>}[d] & X \we Y \ar@{<=>}[d] & \\
		             & \k (X) \se \k (Y) \ar@{=>}[r] & |\k (X)| \we |\k (Y)| \ar@{<=>}[r] & |\k (X)| \he |\k (Y)| } 
\end{displaymath}

\begin{displaymath}
\xymatrix@C=20pt{ \x (K) \he \x (L)  \ar@{=>}[r] & \x (K) \se \x (L) \ar@{=>}[r] \ar@{<=>}[d] & \x (K) \we \x (L) \ar@{<=>}[d] & \\
		             & K \se L \ar@{=>}[r] & |K| \we |L| \ar@{<=>}[r] & |K| \he |L| } 
\end{displaymath}

Here $\he$ denotes that the spaces are homotopy equivalent.

The Wallet $W$ satisfies $W\searrow *$, however $W \overset{he}{\simeq \! \! \! \! \! \! \! \: /} *$. Therefore $X\se Y \Rightarrow \! \! \! \! \! \! \! \! \: / \ \ X\he Y$.

Since $|K|\he |L| \Rightarrow \! \! \! \! \! \! \! \! \: / \ \ K\se L$, we also have that $X \we Y \Rightarrow \! \! \! \! \! \! \! \! \: / \ \ X\se Y$.

\bigskip 

Note that, if $X\we Y$ and their Whitehead group $Wh (\pi_1 (X))$ is trivial, then $|\k (X)|$ and $|\k (Y)|$ are homotopy equivalent CW-complexes with trivial Whitehead group and therefore, simple homotopy equivalent. It follows from \ref{main} that $X\se Y$. Thus we have proved

\begin{coro}
Let $X$, $Y$ be weak equivalent finite $T_0$-spaces such that $Wh (\pi_1 (X))=0$. Then $X\se Y$.
\end{coro}

As another immediate consequence of the theorem, we have

\begin{coro}
Let $X$, $Y$ be finite $T_0$-spaces. If $X\searrow Y$, then $X'\searrow Y'$.
\end{coro}

Note also that, as a corollary of the theorem one deduces the following known fact: if $K$ and $L$ are finite simplicial complexes such that $K\searrow L$, then $K'\searrow L'$.

\bigskip

\bigskip

\textbf{Collapsible finite spaces}

\bigskip

As one can imagine, we will say that a finite $T_0$-space $X$ is \textit{collapsible} if $X\searrow *$.\ 

Observe that every contractible finite $T_0$-space is collapsible, however the converse is not true. The Wallet $W$ introduced in \ref{wallet} is collapsible and non-contractible. 

One could ask if there is a finite $T_0$-space which is homotopically trivial but non-collapsible. We will come back to this question in a minute.

\bigskip

Note that if a finite $T_0$-space $X$ is collapsible, its associated simplicial complex $\k (X)$ is also collapsible. Moreover, if $K$ is a collapsible complex, then $\x(K)$ is a collapsible finite space. Therefore, if $X$ is a collapsible finite space, its subdivision $X'$ is also collapsible.

\bigskip

Let us consider now a compact contractible polyhedron $X$ such that any triangulation of $X$ is non-collapsible, for instance the Dunce Hat \cite{Zee}.

Let $K$ be any triangulation of $X$. Now we claim that $\x (K)$ is homotopically trivial because $X$ is contractible. Nevertheless, $\x (K)$ cannot be collapsible since $K'$ is not collapsible.

\bigskip

This is an interesting example of how classical simple homotopy theory can help us to answer natural questions on finite spaces. 
It was not easy to find a non-contractible homotopically-trivial finite space as we did in \ref{wallet}. However, the finite space associated to the Dunce Hat, despite being much bigger than $W$, arises in a more natural way.

\bigskip

We have the following situation $$\textrm{contractible} \Rightarrow \textrm{collapsible} \Rightarrow \textrm{homotopically trivial}$$ and none of the converses holds.

\bigskip

\bigskip

\textbf{Minimal simple models}

\bigskip

The core of a finite space $X$ is the smallest space which is homotopy equivalent to $X$. As we pointed out before the core is unique up to homeomorphism.

In \cite{Bar} we have studied the \textit{minimal finite models} of a space $X$ (not necessarily finite), which are the smallest spaces which are 
weak equivalent to $X$. We proved that in general these models are not unique. In \cite{Bar} we characterized the minimal finite models of finite graphs (finite
CW-complexes of dimension one) and as an example we showed that $\bigvee\limits_{i=1}^3 S^1$ has three minimal finite models up to homeomorphism.

It has sense to make the following definition.

\begin{defi}
A \textit{minimal simple model} of a finite $T_0$-space $X$ is a finite $T_0$-space simply equivalent to $X$ of minimum cardinal. We will  say that a space is a minimal simple model if it is a minimal simple model of itself.
\end{defi}

For finite $T_0$-spaces we have that $$\textrm{minimal finite model} \Rightarrow \textrm{minimal simple model} \Rightarrow \textrm{minimal finite space}$$

Note that if the Whitehead group $Wh(\pi_1(X))$ is trivial, the converse of the first implication holds. 

Therefore if we have a finite $T_0$-space $X$ such that $Wh(\pi_1(X))=0$, we could reach any minimal finite model of $X$ ``simply'' by adding and removing weak points from $X$.

\smallskip

Going back to the first paragraph of this section: Elementary collapses and expansions give us an effective method of reduction when the space has trivial Whitehead group. Unfortunately this is not exactly what one looks for since it is not possible to get a minimal simple model just by taking away weak points. More explicitly:

A minimal simple model has no weak points. However, if we consider a triangulation $K$ of the Dunce Hat, and we remove as many weak points as we can from $\x (K)$, we will obtain a space without weak points which is not a minimal simple model.

This is very different to the homotopy type case, where removing beat points is an effective way of getting a core.

\bigskip

Of course there is not uniqueness of minimal simple models. For example we can consider the space $\mathbb{S} D_3$

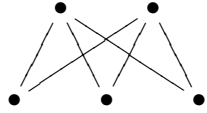
\begin{displaymath}
\xymatrix@C=6pt{ & \bullet \ar@{-}[dl] \ar@{-}[dr] \ar@{-}[drrr] & & \bullet \ar@{-}[dlll] \ar@{-}[dl] \ar@{-}[dr]\\ 
		\bullet & & \bullet & & \bullet} 
\end{displaymath}

and its opposite, which are minimal simple models because they are minimal finite models. Moreover $\mathbb{S} D_3\se (\mathbb{S} D_3) ^{op}$ since $\k (\mathbb{S} D_3)$ and $\k ((\mathbb{S} D_3)^{op})$ are isomorphic.

\bigskip

\section{Simple Homotopy Equivalences: The Second Main Theorem}

In this section we prove the second main result of the article. We present the notion of  a \textit{distinguished} map between finite $T_0$-spaces and the 
class $\dd$ generated by these maps. 

We shall prove that every map in $\overline{\d}$ induces a simple homotopy equivalence between the associated complexes. Conversely, if $|\k(f)|$ 
is a simple homotopy equivalence, then $f\in \overline{\d}$.

Furthermore, for a simplicial map $\varphi:K\to L$, $|\varphi|$ is a simple homotopy equivalence if and only if $\x(\varphi)\in \overline{\d}$.

\bigskip

Recall first that a homotopy equivalence $f:|K|\to |L|$ between compact polyhedra is a simple homotopy equivalence if it is homotopic to a 
composition of a finite sequence of maps $|K| \to |K_1|\to\ldots \to |K_n|\to |L|$, each of them an expansion or a homotopy inverse of one \cite{Coh,Sie}.

\smallskip

We prove first that homotopy equivalences between finite spaces induce simple equivalences between the associated polyhedra.

\begin{teo} \label{ehinduce}
If $f:X\to Y$ is a homotopy equivalence between finite $T_0$-spaces, then $|\k(f)|: |\k (X)|\to |\k (Y)|$ is a simple homotopy equivalence.
\end{teo} 
\begin{proof}
Let $X_c$ and $Y_c$ be cores of $X$ and $Y$. Let $i_X:X_c\to X$ and $i_Y:Y_c\to Y$ be the inclusions and $r_X:X\to X_c$, $r_Y:Y\to Y_c$ retractions of $i_X$ and $i_Y$ such that $i_Xr_X\simeq 1_X$ and $i_Yr_Y\simeq 1_Y$.

Since $r_Yfi_X:X_c\to Y_c$ is a homotopy equivalence between minimal finite spaces, it is a homeomorphism. Therefore $\k(r_Yfi_X):\k (X_c)\to \k(Y_c)$ is an isomorphism and then $|\k(r_Yfi_X)|$ is a simple equivalence.

\bigskip

Since $\k(X)\searrow \k(X_c)$, $|\k(i_X)|$ is a simple equivalence, and then the homotopy inverse $|\k(r_X)|$ is also a simple equivalence. Analogously $|\k(i_Y)|$ is a simple equivalence.

Finally, since $f\simeq i_Yr_Yfi_Xr_X$, it follows that $|\k(f)|\simeq |\k(i_Y)| |\k(r_Yfi_X)| |\k(r_X)|$ is a simple equivalence.
\end{proof}

\begin{defi}
A map $f:X\to Y$ between finite $T_0$-spaces is \textit{distinguished} if $f^{-1}(U_y)$ is contractible for each $y\in Y$. We denote by $\d$ 
the class of distinguished maps.
\end{defi}

\begin{obs} \label{eind}
From the proof of \ref{weak}, it is clear that if $x\in X$ is a weak point such that $U_x\smallsetminus \{x\}$ is contractible, the inclusion $X\smallsetminus \{x\}\hookrightarrow X$ is distinguished.
\end{obs}

Note that by the theorem of McCord (\cite{Mcc}; Theorem 6), every distinguished map is a weak equivalence and therefore induces a homotopy equivalence between the associated complexes. We will prove in \ref{dis} that in fact the induced map is a simple equivalence.

\begin{obsi} \label{hardie}
Let $X$ be a finite $T_0$-space. In \cite{Har} and \cite{May2} it is proved that the map $h:X'\to X$ defined by $h(C)=max(C)$ is a weak equivalence. 
Moreover $\k(h):\k (X)'\to \k(X)$ is a simplicial approximation to the identity and then $\k(h)$ is in fact a simple homotopy equivalence. 
We give here a different approach: Since $$h^{-1}(U_x)=\{C \ | \ max(C)\le x \}=\x (\k (U_x))=\x (x\k (U_x\smallsetminus \{x\}))$$ 
is contractible for every $x\in X$, $h$ is a distinguished map. In particular, by \ref{dis}, $h$ induces a simple homotopy equivalence.
\end{obsi}

It is easy to see that homeomorphisms are distinguished: If $f:X\to Y$ is a homeomorphism, then $f^{-1}(U_y)=U_{f^{-1}(y)}$, which is contractible.

However homotopy equivalences are not distinguished in general. The map

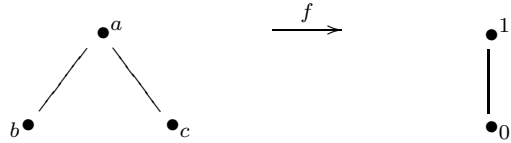
\begin{displaymath}
\xymatrix@C=6pt{ & \ \ \bullet ^a \ar@{-}[dl] \ar@{-}[dr] & \\
                  _b \bullet \ & & \ \bullet _c } 
\qquad
\xymatrix@C=10pt{ \ar@{->}^f[rr] & & \\
									& & & }
\qquad
\xymatrix@C=6pt{ \ \ \bullet ^1 \ar@{-}[d] \\
                 \ \ \bullet _0 }             
\end{displaymath}

defined by $f(a)=1$, $f(b)=f(c)=0$ is a homotopy equivalence because both spaces are contractible, however $f^{-1}(U_0)=\{b,c\}$ which is not contractible. 

\begin{teo} \label{dis}
Every distinguished map induces a simple homotopy equivalence. 
\end{teo}
\begin{proof}
Suppose $f:X\to Y$ is distinguished. We define the \textit{non-Hausdorff mapping cylinder} $B(f)$ as the set $X\sqcup Y$ with the following order. Given $a, b$ in $B(f)$, $a\le b$ if one of the following holds:

\begin{itemize}
\item $a,b\in X$ and $a\le b$ in $X$.
\item $a,b\in Y$ and $a\le b$ in $Y$.
\item $a\in X$, $b\in Y$ and $f(a)\le b$ in $Y$. 
\end{itemize}

Note that the space $B(X)$ constructed in \ref{bedex} is a particular case of the non-Hausdorff cylinder when $f=h:X'\to X$ as defined in \ref{hardie}.

We will show that both $X$ and $Y$ expand to $B(f)$.

\bigskip

Labeling all the elements $x_1, x_2, \ldots, x_n$ of $X$ in such a way that $x_i\le x_j$ implies $i\le j$ and defining $Y_i=\{x_1, x_2, \ldots, x_i\} \sqcup Y \subseteq B(f)$ for every $0\le i\le n$, we have that $$\overline{\{x_i\}}^{Y_i}\smallsetminus \{x_i\}=\{y\in Y \ | \ y\ge f(x_i)\}$$

is homeomorphic to the contractible space $\overline{\{f(x_i)\}}^Y$. Therefore $Y_i \ce Y_{i-1}$ for $1\le i\le n$, and then $Y=Y_0$ expands to $B(f)=Y_n$.

Notice that we have not used yet the fact that $f$ is distinguished.

\bigskip

Now order the elements $y_1,y_2, \ldots, y_m$ of $Y$ in such a way that $y_i\le y_j$ implies $i\le j$ and define $X_i=X \sqcup \{y_{i+1},y_{i+2}, \ldots, y_m\} \subseteq B(f)$ for every $0\le i\le m$. Then $$U_{y_{i}}^{X_{i-1}}\smallsetminus \{y_{i}\}=\{x\in X \ | \ f(x)\le y_{i}\},$$ is homeomorphic to $f^{-1}(U_{y_i})$, which is contractible by hypothesis. Thus $X_{i-1}\ce X_i$ for $1\le i\le m$. Therefore $B(f)=X_0$ collapses to $X=X_m$.

\bigskip

The following diagram
\begin{displaymath}
\xymatrix@C=30pt{ & B(f) & \\
                 X \ar@{^{(}->}^{i_X}[ur] \ar@{->}^f[rr] & & Y \ar@{_{(}->}_{i_Y}[ul] }             
\end{displaymath}
does not commute, but $i_X\le i_Yf$ and then $i_X\simeq i_Yf$. Therefore $|\k(i_X)|\simeq |\k(i_Y)| |\k(f)|$. The expansions $|\k(i_X)|$ and $|\k(i_Y)|$ are simple equivalences, and then so is $|\k(f)|$.
\end{proof}

We have already showed that expansions, homotopy equivalences and distinguished maps induce simple equivalences at the level of complexes. 

Note that if $f,g,h$ are three maps between finite $T_0$-spaces such that $fg\simeq h$ and two of them induce simple equivalences, then the third map also does.

This provides us a way to construct new maps that induce simple equivalences:

\begin{defi}
Let $\mathcal{C}$ be a class of continuous maps between finite $T_0$-spaces. We define recursively the class $\mathcal{C}_n$ in the following way. $\mathcal{C}_0=\mathcal{C}$, $$\mathcal{C}_{n+1}=\{f,g,h \ | \ fg\simeq h \textrm{ and such that 2 of the 3 are in }\mathcal{C}_n  \}$$

We call $\overline{\mathcal{C}}=\bigcup\limits_{n\in \mathbb{N}}\mathcal{C}_n$, the class \textit{generated} by $\mathcal{C}$.  

In other words, $\overline{\mathcal{C}}$ is the smallest class having the 2-of-3 property (up to homotopy) containing $\mathcal{C}$.  
\end{defi}

It is clear that no matter what $\c$ is, the class $\overline{\c}$ is closed by composition and homotopy.

Note that if every map in $\c$ induces a simple equivalence, each map in $\overline{\c}$ does. If we denote $\e$ the class of elementary 
expansions of finite $T_0$-spaces, it is clear then that every map in $\overline{\e}$ induces a simple homotopy equivalence.

\bigskip

Observe that the class of simple equivalences between CW complexes is the smallest class closed by the 2-of-3 property (up to homotopy) 
containing elementary expansions. In the setting of finite spaces, a map which induces a simple equivalence needs not have a homotopy inverse. 
This is the reason why the definition of $\overline{\e}$ is not as simple as in the setting of CW-complexes.

\begin{obs} \label{dine}
Every expansion of finite $T_0$-spaces is in $\overline{\e}$ because it is a composition of maps in $\e$. The proof of \ref{dis} shows that for every distinguished map $f$, there exist expansions $i$, $j$ such that $i\simeq jf$. Therefore $f\in \overline{\e}$.
\end{obs}

A map $f:X\to Y$ such that $f^{-1}(\overline{\{y\}})$ is contractible for every $y$ needs not be distinguished. However it is clear that $|\k(f)|$ is a simple equivalence since $f^{op}:X^{op}\to Y^{op}$ is distinguished. Here, $f^{op}$ denotes the map that coincides with $f$ in the underlying sets. We prove that in fact $f\in \overline{\d}$.

\medskip

We denote $\d ^{op}$ the class of maps $f$ such that $f^{op}\in \d$.

\begin{prop} \label{dopind}
Let $f:X\to Y$ be in $\d ^{op}$. Then $f\in \overline{D}$.
\end{prop}
\begin{proof}
Consider the following commutative diagram
\begin{displaymath}
\xymatrix@C=30pt{ X \ar@{->}^f[r] & Y \\
 X'=(X^{op})' \ar@{->}^{h_X}[u] \ar@{->}^{f'}[r] \ar@{->}^{h_{X^{op}}}[d] & Y'=(Y^{op})' \ar@{->}^{h_Y}[u] \ar@{->}^{h_{Y^{op}}}[d] \\
                  X^{op} \ar@{->}^{f^{op}}[r] & Y^{op} }             
\end{displaymath}
Here, $f'$ denotes the map $\x(\k(f))$. Since $\overline{\d}$ satisfies the 2-of-3 property and $h_{X^{op}}$, $h_{Y^{op}}$, $f^{op}$ are distinguished by \ref{hardie}, then $f'\in \overline{\d}$. And since $h_X$, $h_Y$ are distinguished, $f\in \overline{\d}$.
\end{proof}

Therefore if $x\in X$ is a weak point such that $\overline{\{x\}}\smallsetminus \{x\}$ is contractible, the inclusion $X\smallsetminus \{x\} \hookrightarrow X$ lies in $\overline{\d}$. 

\begin{coro} \label{todoigual}
By \ref{eind} and what we have just proved, every elementary expansion is in $\overline{\d}$, which proves that $\overline{\e}\subseteq \overline{\d}$. 
From \ref{dine},  it follows that $\overline{\e}= \overline{\d}$.

From \ref{dopind} and a simetrical result we also obtain that $\overline{\d}=\overline{\d ^{op}}$.\

The proof of \ref{ehinduce} shows that if $f:X\to Y$ is a homotopy equivalence, then $fi_X\simeq i_Yr_Yfi_X$ where $i_X$, $i_Y$ are expansions and $r_Yfi_X$ is a homeomorphism. Thus, every homotopy equivalence lies in $\overline{\e}=\overline{\d}=\overline{\d ^{op}}$.
\end{coro}

Below we shall prove that $\overline{\e}=\overline{\d}=\overline{\d ^{op}}$ is exactly the class of maps which induce simple homotopy equivalences between the associated polyhedra.

\begin{lema} \label{lema1}
Let $\varphi,\psi:K\to L$ be simplicial maps which lie in the same contiguity class. Then $\x (\varphi)\simeq \x (\psi)$.
\end{lema}
\begin{proof}
Assume that $\varphi$ and $\psi$ are contiguous. Then the map $f :\x(K)\to \x(L)$, defined by $f (S)=\varphi(S)\cup \psi(S)$ is well-defined and continuous. Moreover $\x(\varphi)\le f \ge \x(\psi)$, and then $\x(\varphi)\simeq \x(\psi)$.
\end{proof}

Given $n\in\mathbb{N}$ we denote by $K^n$ the n-th barycentric subdivision of $K$.

\begin{lema} \label{lema2}
Let $\lambda :K^n\to K$ be a simplicial approximation to the identity. Then $\x(\lambda)\in \overline{D}$.
\end{lema}
\begin{proof}
Suppose first that $n=1$. Let $\lambda :K'\to K$ be any simplicial approximation to $1_{|K|}$. Then $\x(\lambda): \x(K)'\to \x(K)$ is homotopic to $h_{\x(K)}$, for if $S_1\subsetneq S_2\subsetneq \ldots \subsetneq S_m$ is a chain of simplices of $K$, then $\x(\lambda)(\{S_1,S_2,\ldots,S_m\})=\{\lambda(S_1),\lambda(S_2),\ldots,\lambda(S_m)\}\subseteq S_m=h_{\x(K)}(\{S_1,S_2,\ldots,S_m\})$.

By \ref{hardie}, it follows that $\x(\lambda)\in \overline{\d}$.

\bigskip

Now suppose that $n\ge 1$. A composition of approximations to the identity $K^{i+1}\to K^i$ for $0\le i<n$ defines an approximation to the identity $\nu :K^n\to K$. Since $\overline{\d}$ is closed by compositions, $\x(\nu)\in \overline{\d}$.

Any other approximation $\lambda: K^n\to K$ to $1_{|K|}$, is contiguous to $\nu$. By \ref{lema1}, $\x(\lambda)$ is homotopic to $\x(\nu)$ and lies in $\overline{\d}$.

\end{proof}

\begin{lema} \label{lema3}
Let $\varphi,\psi:K\to L$ be simplicial maps such that $|\varphi|\simeq|\psi|$ and such that $\x(\varphi)\in \overline{\d}$. Then $\x(\psi)$ also lies in $\overline {\d}$.
\end{lema}
\begin{proof}
There exist $n\ge 1$ and simplicial approximations $\widetilde{\varphi}, \widetilde{\psi}: K^n\to L$ to $|\varphi|$ and $|\psi|$ in the same contiguity class. 

Let $\lambda :K^n\to K$ be a simplicial approximation to $1_{|K|}$. Then $\varphi\lambda$ is an approximation to $|\varphi|$ and therefore $\varphi\lambda$ and $\widetilde{\varphi}$ are contiguous. Analogously $\psi\lambda$ and $\widetilde{\psi}$ are contiguous. Hence, $\varphi\lambda$ and $\psi\lambda$ lie in the same contiguity class and, by \ref{lema1}, $\x(\varphi)\x(\lambda)= \x(\varphi\lambda)\simeq \x(\psi\lambda) =\x(\psi)\x(\lambda)$. By \ref{lema2} $\x(\varphi),\x(\lambda)\in \overline{\d}$. It follows that $\x(\psi)\in \overline{\d}$.
\end{proof}

\begin{teo}
Let $K_0,K_1,\ldots,K_n$ be finite simplicial complexes and let $f_i:|K_i|\to |K_{i+1}|$ be such that for each $0\le i<n$ one of the following holds:

\begin{enumerate}
	\item[(1)] $f_i=|\varphi_i|$ where $\varphi_i:K_i\to K_{i+1}$ is a simplicial map such that $\x(\varphi_i)\in \overline{\d}$.
	\item[(2)] $f_i$ is a homotopy inverse of $|\varphi_i|$ where $\varphi_i:K_{i+1}\to K_{i}$ is a simplicial map such that $\x(\varphi_i)\in \overline{\d}$.
\end{enumerate}
Let $\varphi:K_0\to K_n$ be a simplicial map such that $|\varphi|\simeq f_{n-1}f_{n-2}\ldots f_0$. Then $\x(\varphi)\in \overline{\d}$.
\end{teo} 
\begin{proof}
We may assume that $f_0$ satisfies condition $(1)$. Otherwise we define $\widetilde{K_0}=K_0$, $\widetilde{f_0}=|1_{K_0}|:|\widetilde{K_0}|\to |K_0|$ and then $|\varphi|\simeq f_{n-1}f_{n-2}\ldots f_0\widetilde{f_0}$.

\bigskip

The theorem is proved by induction on $n$.

If $n=1$, $|\varphi|\simeq |\varphi_0|$ where $\x(\varphi_0)\in \overline{\d}$ and the result follows straightforward from \ref{lema3}.

Suppose now that the theorem is true for $n$. Let $K_0,K_1,\ldots, K_n,K_{n+1}$ be finite complexes, $f_i:|K_i|\to |K_{i+1}|$ maps satisfying conditions $(1)$ or $(2)$ and $\varphi:K_0\to K_{n+1}$ such that $|\varphi|\simeq f_nf_{n-1}\ldots f_0$.

We consider two cases: $f_n$ satisfies condition $(1)$ or $f_n$ satisfies condition $(2)$.

\bigskip

In the first case we define $g:|K_0|\to |K_n|$, $g=f_{n-1}f_{n-2}\ldots f_0$.

Let $\widetilde{g}:K_0^m\to K_n$ be a simplicial approximation to $g$ and let $\lambda :K_0^m\to K_0$ be a simplicial approximation to the identity. Then $|\widetilde{g}|\simeq g|\lambda|=f_{n-1}f_{n-2}\ldots f_1(f_0|\lambda|)$ where $f_0|\lambda|=|\varphi_0\lambda|$ and $\x(\varphi_0\lambda)=\x(\varphi_0)\x(\lambda)\in \overline{\d}$ by \ref{lema2}. 

By induction, $\x(\widetilde{g})\in \overline{\d}$, and then $\x(\varphi_n\widetilde{g})\in \overline{\d}$. Since $|\varphi\lambda|\simeq f_ng|\lambda|\simeq f_n|\widetilde{g}|=|\varphi_n\widetilde{g}|$, by lemma \ref{lema3}, $\x(\varphi\lambda)$ lies in $\overline{\d}$. Therefore $\x(\varphi)\in \overline{\d}$.

\bigskip

In the other case, $|\varphi_n\varphi|\simeq f_{n-1}f_{n-2}\ldots f_0$ and by induction, $\x(\varphi_n\varphi)\in \overline{\d}$. Therefore $\x(\varphi)$ also lies in $\overline{\d}$.
\end{proof}

\begin{coro}
Let $\varphi:K\to L$ be a simplicial map such that $|\varphi|$ is a simple homotopy equivalence. Then $\x(\varphi)\in \overline{\d}$.
\end{coro}
\begin{proof}
Since $|\varphi|$ is a simple equivalence, there exist finite complexes $K=K_0, K_1,\ldots K_n=L$ and maps $f_i:|K_i|\to |K_{i+1}|$, which are simplicial expansions or homotopic inverses of simplicial expansions, and such that $|\varphi|\simeq f_{n-1}f_{n-2}\ldots f_0$.

By our First Main Theorem \ref{main}, simplicial expansions between complexes induce expansions between the associated finite spaces which lie in $\dd$ by \ref{todoigual}. Therefore, the last theorem applies. 
\end{proof}

Now we are ready to prove the second important result of this article which is the analogue for maps of our First Main Theorem. In fact, we have already done most of the work.

\begin{smteo}
\begin{enumerate}
\item[]
\item[(a)] Let $f:X\to Y$ be a map between finite $T_0$-spaces. Then $f\in \dd$ if and only if $|\k(f)|:|\k(X)|\to |\k(Y)|$ is a simple homotopy equivalence. 

\item[(b)] Let $\varphi:K\to L$ be a simplicial map between finite simplicial complexes. Then $|\varphi|$ is a simple homotopy equivalence if and only if $\x(\varphi)\in \dd$.
\end{enumerate}
\end{smteo} 

\begin{proof}
Suppose $f:X\to Y$ is a map such that $|\k(f)|$ is a simple equivalence. By the last corollary, $f':X'\to Y'$ lies in $\dd$ and since $fh_X=h_Yf'$, we have that $f\in \dd$.

If $\varphi:K\to L$ is a simplicial map such that $\x(\varphi)\in \dd$, then $|\varphi'|:|K'|\to |L'|$ is a simple equivalence. Here $\varphi'=\k(\x(\varphi))$ is the barycentric subdivision of $\varphi$.

Let $\lambda_K:K'\to K$ and $\lambda_L:L'\to L$ be simplicial approximations to the identities. Then $\lambda_L\varphi'$ and $\varphi\lambda_K$ are contiguous. In particular $|\lambda_L||\varphi'|\simeq |\varphi||\lambda_K|$ and then $|\varphi|$ is a simple equivalence.
\end{proof}

Since there are homotopy equivalences which are not simple in the setting of polyhedra, the theorem says that in the setting of finite spaces the inclusions $$\{homotopy \ equivalences\}\subsetneq \dd \subsetneq \{weak \ equivalences\}$$ are both strict. However, if $f:X\to Y$ is a weak equivalence between finite $T_0$-spaces with trivial Whitehead group, then $f\in \dd$.

\email{jbarmak@dm.uba.ar, gminian@dm.uba.ar}


\begin{thebibliography}{99}

\bibitem{Ale} P.S. Alexandroff. \textit{Diskrete R\"aume}.
    MathematiceskiiSbornik (N.S.) 2(1937), 501-518.

\bibitem{Bar} J.A. Barmak and E.G. Minian. \textit{Minimal finite models}.
    Preprint (2006).

\bibitem{Coh} M.M. Cohen. \textit{A Course in Simple Homotopy Theory}.
    Springer-Verlag New York, Heidelberg, Berlin (1970).

\bibitem{Har} K.A. Hardie and J.J.C. Vermeulen. \textit{Homotopy theory of finite and locally finite $T_0$-spaces}.
    Exposition Math. 11 (1993), 331-341.

\bibitem{May} J.P. May. \textit{Finite topological spaces}.
    Notes for REU (2003).

\bibitem{May2} J.P. May. \textit{Finite spaces and simplicial complexes}.
    Notes for REU (2003).

\bibitem{Mcc} M.C. McCord. \textit{Singular homology groups and homotopy groups of finite topological spaces}.
    Duke Mathematical Journal 33(1966), 465-474.

\bibitem{Mil} J. Milnor. \textit{Whitehead Torsion}. Bull. AMS 72(1966), 358-426.

\bibitem{Osa} T. Osaki. \textit{Reduction of finite topological spaces}.
    Interdiciplinary Information Sciences 5(1999), 149-155.

\bibitem{Sie} L.C. Siebenmann. \textit{Infinite simple homotopy types}. Indagationes Math. 32(1970), 479-495.

\bibitem{Sto} R.E. Stong. \textit{Finite topological spaces}.
    Trans. Amer. Math. Soc. 123(1966), 325-340.

\bibitem{Whi} J.H.C Whitehead. \textit{Simplicial spaces, nuclei and m-groups}.
    Proc. London Math. Soc. 45(1939), 243-327.

\bibitem{Whi2} J.H.C Whitehead. \textit{On incidence matrices, nuclei and homotopy types}.
    Ann. of Math. 42(1941), 1197-1239.

\bibitem{Whi4} J.H.C Whitehead. \textit{Simple homotopy types}.
    Amer. J. Math. 72(1950), 1-57.

\bibitem{Zee} E.C. Zeeman \textit{On the dunce hat}.
    Topology 2(1964), 341-358.

\end{thebibliography}
\end{document}